\documentclass[a4paper,12pt,reqno]{amsart}

\usepackage{lscape}
\usepackage{amsmath,amsfonts,amsthm,amssymb, latexsym,euscript,dsfont,graphicx,units,color}
\usepackage[latin1]{inputenc}
\usepackage[font=footnotesize]{subfig}
\usepackage{pdfsync}

\newcommand{\1}{1 \hspace*{-0.2ex}\rule{0.10ex}{1.5ex}\hspace{0.2ex}}

\newcommand{\E}{\mathbb E}

\newcommand{\R}{\mathbb R}

\newcommand{\bean}{\begin{eqnarray*}}
\newcommand{\eean}{\end{eqnarray*}}
\newcommand{\ben}{\begin{enumerate}}
\newcommand{\een}{\end{enumerate}}
\newcommand{\beq}{\begin{equation}}
\newcommand{\eeq}{\end{equation}}

\newtheorem{theorem}{Theorem}[section]

\newtheorem{corollary}[]{Corollary}

\newtheorem{lemma}[theorem]{Lemma}

\newtheorem{proposition}[]{Proposition}

\theoremstyle{remark}
\newtheorem{remark}[]{Remark}

\begin{document}
\title[]{\bf Certainty bands for the conditional cumulative distribution function and applications}

\author[Ferrigno, Foliguet, Maumy-Bertrand Muller-Gueudin]{By S. Ferrigno$^1$\and B. Foliguet$^2$\and M. Maumy-Bertrand$^3$\and A. Muller-Gueudin$^1$}

\begin{abstract}
In this paper, we establish uniform asymptotic certainty bands for the conditional cumulative distribution function. 
To this aim, we give exact rate of strong uniform consistency for the local linear estimator of this function. 
The corollaries of this result are the asymptotic certainty bands for the quantiles and the regression function. We illustrate our results with simulations and an application on fetopathologic data.
\end{abstract}
\maketitle

\small
\noindent$^{1}$IECL (UMR 7502), Nancy-Université, CNRS, INRIA, members of the BIGS (BIology, Genetics and Statistics) team at INRIA, France\footnote{{\tt{sandie.ferrigno@univ-lorraine.fr}},{\tt{aurelie.gueudin@univ-lorraine.fr}}}, $^{2}$Service de foetopathologie et de placentologie, Maternité Régionale Universitaire, CHU Nancy, France\footnote{{\tt{b.foliguet@maternite.chu-nancy.fr}}}, $^{3}$IRMA (UMR 7501), Université de Strasbourg, France.\footnote{ {\tt {mmaumy@math.unistra.fr}}}\\

\subsection*{Keywords:} Conditional cumulative distribution function, local polynomial estimator, uniform asymptotic certainty bands, regression function, quantiles.
 
\section{\bf Introduction}\label{Intro}
\subsection{Motivations}

Consider $(X,Y)$, a random vector defined in $\mathbb{R}\times \mathbb{R}$. 
Here $Y$ is the variable of interest and $X$ the concomitant variable. 
Throughout, we work with a sample $\{(X_i,Y_i)_{1\leqslant i \leqslant n}\}$ of independent and identically replica of $(X,Y)$. 
We will assume that $(X,Y)$ [resp. $X$] has a density function $f_{X,Y}$ [resp. $f_{X}$] with respect to the Lebesgue measure. 
In this paper, we will mostly focus on the regression function of $\psi({ Y})$ evaluate at $X=x$ defined by:
\begin{equation}\label{E:1}
m_{\psi}({ x})=\mathbb{E}\left(\psi({Y})|{ X}={ x}\right)=\frac{1}{f_{ X}({ x})}\int_{\mathbb{R}}\psi({y})f_{X,Y}({ x},{ y})d{ y},\quad\mbox{with  }f_X(x)\neq 0
\end{equation}
\noindent whenever this regression function is meaningful. 
Here and elsewhere, $\psi$ denotes a specified measurable function, which is assumed to be bounded on each compact subinterval of $\mathbb{R}$.
\vspace{0.3cm}

Because of numerous applications, the problem of estimating the function $m_{\psi}$, the density function $f_{X}$ and the regression function $m_{\mathbb{\psi={\mbox{Id}}}}$
has been the subject of considerable interest during the last decades. We can cite for example Nadaraya~\cite{Nadaraya}, Watson~\cite{Watson}, Devroye~\cite{Devroye}, Collomb~\cite{Collomb}, H\"ardle~\cite{Hardle} and specially mention two articles, Einmahl and Mason~\cite{Einmahl} and Deheuvels and Mason~\cite{DeheuvelsMason} for two reasons. 
The first is that these articles study an estimator of $m_{\psi}$ and its properties. 
The second is that we use the tools which are developed in these articles in order to establish our proofs.
We now choose $\psi=\psi_t$ defined by $\psi_t(y)=\1_{\{{ y}\leqslant { t}\}}\mbox{ with }{ t}\in\mathbb{R}$ arbitrary but fixed,
and $\1$ the indicator function, so we obtain the conditional cumulative distribution function ({\it cond-cdf}) of $Y$ given $X=x$, defined by: 
\begin{equation}\label{E:3}
\forall t\in \mathbb{R},\quad F({t}|{x})=m_{\psi_t}(x)=\mathbb{E}\left(\1_{\{{ Y}\leqslant { t}\}}|{ X}={ x}\right)
=\mathbb{P}\left({Y}\leqslant {t}|{X}={x}\right).
\end{equation}
Saying that, we are implicitly assuming the existence of a regular version
for the conditional distribution of $Y$ given $X$. 

\vspace{0.3cm}
In this article, we study the conditional cumulative distribution function and a nonparametric estimator associated to this function. 
The {\it cond-cdf} has the advantages of completely characterizing the law of the random considered variable, allowing to obtain the regression function, the density function, the moments and the conditional quantile function. The {\it cond-cdf} is also used for example, in medicine (Gannoun {\emph {et al.}}~\cite{Gannoun}) or econometric domain (Li {\emph {et al.}}~\cite{Li}).
\vspace{0.3cm}

Introduce the Nadaraya-Watson estimator (see Nadaraya~\cite{Nadaraya} and Watson~\cite{Watson}) of the {\it cond-cdf} $F({t}|{x})$, for all ${t}\in \mathbb{R}$ and ${x}\in \mathbb{R}$, defined by: 
\begin{equation}\label{E:4}
\widehat{F}^{(0)}_{n}({t,h_n}|{x})=\frac{\sum_{i=1}^n \1_{\{{Y}_i\leqslant { t}\}} { K}\left(\frac{{ x}-{{X}}_i}{h_n}\right)}{\sum_{i=1}^n  { K}\left(\frac{{ x}-{{X}}_i}{h_n}\right)} \textrm{ for }  \sum_{i=1}^n  { K}\left(\frac{{ x}-{{X}}_i}{h_n}\right)\neq0
\end{equation}
where ${ K}(\cdot)$ is a positive-valued kernel function defined on $\mathbb{R}$ 
and $(h_n)_{n\geqslant 1}$ is the bandwidth, and denotes a non-random sequence of positive constants satisfying some assumptions which will be defined latter.\\ 
For the study of the convergence rate of this estimator, it will be convenient to center $\widehat{F}^{(0)}_{n}({t,h_n}|{x})$ by the estimator of 
$\mathbb{E}(\widehat{F}^{(0)}_{n}({t,h_n}|{x}))$:
\begin{equation}
\widehat{\mathbb{E}}\left(\widehat{F}^{(0)}_{n}({t,h_n}|{x})\right)=\frac{\mathbb{E}\left(\1_{\{Y\leqslant t\}}K\left(\frac{x-X}{h_n}\right)\right)}
{\mathbb{E}\left(K\left(\frac{x-X}{h_n}\right)\right)}\cdot
\end{equation}

\vspace{0.3cm}

\begin{remark}
\hspace{0.1cm}
\begin{enumerate}
\item With the hypothesis we have in mind (see in Section \ref{Hypotheses}), the denominator of this quantity does not cancel. 
\item In general, $\widehat{\mathbb{E}}(\widehat{F}^{(0)}_{n}({t,h_n}|{x}))$ does not coincide with  $\mathbb{E}(\widehat{F}^{(0)}_{n}({t,h_n}|{x}))$. 
However, under mild regularity assumptions, the difference between these two (non-random) quantities becomes asymptotically
negligible as $h_n\searrow 0$ together with $nh_n\nearrow +\infty$ as $n\rightarrow+\infty$.
\end{enumerate}
\end{remark}
\vspace{0.3cm}

The estimator $\widehat{F}^{(0)}_{n}({t,h_n}|{x})$ of the {\it cond-cdf}  has been first treated by Collomb~\cite{Collomb80}.
He proved consistency results, without rates, which are uniform in $x$ and pointwise in $t$.
A Glivenko-Cantelli type theorem for the Nadaraya-Watson estimator, uniform in $t$ and pointwise in $x$, is given in Stute~\cite{Stute}.
Moreover, Stute~\cite{Stute82} was the first to obtain the exact rate of strong uniform consistency on compact intervals for the kernel density estimator. 
The best results for the other estimators in terms of approximate rates that we are aware of are due to H\"ardle {\emph {et al.}}~\cite{Hardle88}.
Such a result plays a fundamental role in obtaining strong uniform  consistency rates in other statistics problems, 
like in the conditional quantile function.
\vspace{0.3cm}

More later, in 2000, Einmahl and Mason~\cite{Einmahl} have determined, under mild regularity conditions on the joint and marginal density functions and under hypotheses 
on the bandwidth $(h_n)$, exact rates of strong uniform  consistency 
for the {\it cond-cdf}. We recall here their result:
\vspace{0.3cm}

\begin{corollary}\label{C:1}(see Corollary 2 in~\cite{Einmahl}.)
Let $I$ be a compact interval.
Assume that $f_{X,Y}$ and $f_X$ satisfy some regularity conditions 
and moreover that $h_n$ satisfies $h_n\searrow 0$, $nh_n\nearrow+\infty, \log h_n^{-1}/\log\log n\rightarrow+\infty$ 
and $nh_n/\log n\rightarrow+\infty$ as $n\rightarrow+\infty$. 
Then we have for any kernel $K$ defined in~\cite{Einmahl}, with probability 1:
\begin{equation}\label{E:lim}
\lim_{n\rightarrow+\infty}\sup_{t\in\mathbb{R}}\sup_{x\in I}\sqrt{\frac{nh_n}{\log (h_n^{-1})}}\left|\widehat{F}^{(0)}_{n}(t,h_n|x)-\widehat{\mathbb{E}}\left(\widehat{F}_{n}^{(0)}(t,h_n|{ x})\right)\right|=
\displaystyle\frac{||{K}||_2}{\displaystyle\sqrt{2\inf_{{x}\in{ I}}f_{ X}({ x})}}
\end{equation}
where
$||K||_2^2=\int_{\mathbb{R}}K^2(u)du$. 
\end{corollary}
\begin{remark}
\hspace{0.1cm}
\begin{enumerate}
\item
Under the assumptions of this corollary, 
the limit in Equation \eqref{E:lim} does not depend on the distribution of the random variable $Y$. 
\item
In 2004, Blondin establishes in~\cite{Blondin} a similar result of the Corollary~\ref{C:1} in the multivariate case, {\emph{i.e.}} $(X,Y)$ is in $\mathbb{R}^r\times\mathbb{R}^d$, $r,d\in\mathbb N ^*$.
\end{enumerate}
\end{remark}
\vspace{0.3cm}

In 2005, Einmahl and Mason~\cite{Einmahl05} have given an uniform in bandwidth consistency of kernel-type function estimators, in the case where $({X},Y)$ is in $\mathbb{R}^r\times\mathbb{R}$, $r\in\mathbb N^*$, and specially for the estimator $\widehat{F}_{n}^{(0)}({t,h_n}|{ x})$ defined in \eqref{E:4}.
We recall below their result:
\begin{theorem}(see Theorem 3 in~\cite{Einmahl05}.)
Let $I$ be a compact subset of $\R^r$ and let $K$ be a kernel defined in~\cite{Einmahl05}. 
Suppose that $f_X$ is continuous and strictly positive on $J$, which is a compact subset of $\R^r$ and contains $I$.
Then, with probability 1, we have for large enough $c>0$ and any $b_n\searrow0$:
\begin{equation}
\limsup_{n\rightarrow+\infty}\sup_{c\log n\leqslant h\leqslant b_n}\sup_{t\in\mathbb{R}}\sup_{x\in I}\sqrt{\frac{nh}{\log(h^{-1})\vee \log\log n}}
\left|\widehat{F}_{n}^{(0)}(t,h|x)-\widehat{\E}\left(\widehat{F}_{n}(t,h|x)\right)\right|<+\infty.
\end{equation}
\end{theorem}
\begin{remark}
In this result, the exact value of the limit is unknown.
\end{remark}
\vspace{0.5cm}

It is a well-known fact the asymptotic bias of the Nadaraya-Watson estimator has a bad form.
To overcome this problem, there exists an alternative: the local polynomial techniques described in Fan and Gijbels~\cite{Fan} or in Tsybakov~\cite{Tsybakov}.

To study the local polynomial estimators, either we can use the $U$-statistics, for example Mint El Mouvid~\cite{Mint El Mouvid}.
But this method implies heavy calculations. Or we can use the empirical processes, for example Dony {\emph{et al.}}~\cite{DonyEinmahlMason}. 
But the results on the empirical processes indexed by classes of functions are established only for classes of real-valued functions.
\vspace{0.3cm}

The present paper is organized as follows. First, we introduce the   local linear estimator of the {\it cond-cdf}, with the main notations and assumptions needed for our task. 
Then we establish an uniform law of the logarithm for the local linear estimator of the {\it cond-cdf}  in Section \ref{Sect:Main_Result}. 
In Section \ref{Sect:3}, we show that limit laws of the logarithm are useful in the construction of uniform asymptotic certainty bands for the {\it cond-cdf}, the regression function and the conditional quantile function. 
Such certainty bands are obtained from simulations in Section \ref{appli} and from fetopathologic data in Section \ref{Sect:feto}. 
Finally, Section \ref{proofs} is devoted to the proofs of our results.  

\subsection{Notations and assumptions}\label{Hypotheses}

Let $(X_1,Y_1),(X_2,Y_2),\dots,$ be independent and identically distributed replica of $(X,Y)$ in $\mathbb{R}\times\mathbb{R}$.
Let ${I}=[a,b],{J}=[a',b']\supsetneq {I}$, two fixed compacts of $\mathbb{R}$.
\vspace{0.3cm}

\noindent 
First, we impose the following set of assumptions upon the distribution of $(X,Y)$:
\begin{enumerate}
\item[(F.1)] $f_{X,Y}$ is continuous on $J\times\mathbb{R}$ and $f_{X}$ is continuous and strictly positive on $J$;
\item[(F.2)] $Y\1_{\{X\in{J}\}}$ is bounded on $\mathbb{R}$.
\end{enumerate}
\begin{remark}
\hspace{0.1cm}
\begin{enumerate}
\item
Under (F.1-2), the {\it cond-cdf} is well defined.
\item
 The assumption (F.2) is very useful for the proof of our results. This boundedness  assumption is common in non-parametric estimation. It ensures the existence of several moments of the {\it cond-cdf}. 
\end{enumerate}
\end{remark}
\vspace{0.3cm}

\noindent 
$K$ denotes a positive-valued kernel function defined on $\mathbb{R}$, fulfilling the conditions:
\begin{enumerate}
\item[(K.1)] ${ K}$ is right-continuous function with bounded variation on $\mathbb{R}$;
\item[(K.2)] ${ K}$ is compactly supported and $\int_{\mathbb{R}}{ K}({ u})\mbox{d}{ u}=1$;
\item[(K.3)] $\int_{\mathbb{R}}{ u}{ K}({ u})\mbox{d}{ u}=0$ and $\int_{\mathbb{R}}{ u}^2{ K}({ u})\mbox{d}{ u}\neq 0$.
\end{enumerate}
We note: $||K||_2^2=\int_{\mathbb{R}}K^2(u)du$.
\vspace{0.3cm}

\noindent 
Further, introduce the following  assumptions on the non-random sequence $(h_n)_{n \geqslant1}$:
\begin{enumerate}
\item[(H.0)] for all $n$, $0<h_n<1$;
\item[(H.1)] $h_n\to 0$, as $n\to +\infty$;
\item[(H.2)] $nh_n/\log n\to+\infty$, as $n\to+\infty$;
\item[(H.3)] $h_n\searrow 0$ and $nh_n\nearrow+\infty$, as $n\rightarrow+\infty$;
\item[(H.4)] $\log(h_n^{-1})/\log\log n\rightarrow+\infty$, as $n\rightarrow+\infty$.
\end{enumerate}

\begin{remark}
\hspace{0.1cm}
\begin{enumerate}
\item The assumption (H.0) is necessary to define $\sqrt{\log(h_n^{-1})}^{-1}$ (see later in our Theorem \ref{theo1}).
\item The assumptions  (H.0-2) are necessary and sufficient for our uniform convergence in probability (see Theorem \ref{theo1}).
\item In order to have almost surely convergence results, we need the assumptions (H.3-4) (see Blondin~\cite{Blondin}).
\item The assumptions (H.0, H.2-4) are called the Csörgö-Révész-Stute assumptions.
\end{enumerate}
\end{remark}
\vspace{0.3cm}

Our aim will be to establish the strong uniform consistency of the local linear estimator of the conditional cumulative distribution function, defined by: 
\begin{equation}
\widehat{F}^{(1)}_{n}(t,h_n|{x})=\frac{\widehat{f}_{n,2}({x,h_n})\widehat{r}_{n,0}(x,t,h_n)-\widehat{f}_{n,1}({x,h_n})\widehat{r}_{n,1}(x,t,h_n)}{\widehat{f}_{n,0}({ x,h_n})\widehat{f}_{n,2}({x,h_n})-\left(\widehat{f}_{n,1}({x,h_n})\right)^2}
\end{equation}
where $^{(1)}$ denotes the order 1 of the local polynomial estimator, and 
\begin{equation}\label{defFnj}
\widehat{f}_{n,j}({x,h_n})=\frac{1}{nh_n}\sum_{i=1}^n\left(\frac{{x}-{ X}_i}{h_n}\right)^j{ K}\left(\frac{{ x}-{ X}_i}{h_n}\right), \mbox{ for }j=0,1,2,
\end{equation}
\begin{equation}\label{defRnj}
\widehat{r}_{n,j}(x,t,h_n)=\frac{1}{nh_n}\sum_{i=1}^n\1_{\{Y_i\leqslant t\}}\left(\frac{{ x}-{ X}_i}{h_n}\right)^j{ K}\left(\frac{{ x}-{ X}_i}{h_n}\right), \mbox{ for }j=0,1.
\end{equation}
\vspace{0.3cm}

\begin{remark}
\hspace{0.1cm}
\begin{enumerate}
\item 
The Nadaraya-Watson estimator $\widehat{F}^{(0)}_{n}({t,h_n}|{x})$ can be also written with the functions $\widehat f_{n,j}$ and $\widehat r_{n,j}$ as 
$$
\widehat{F}^{(0)}_{n}({t,h_n}|{x})=\displaystyle\frac{\widehat{r}_{n,0}({ x},t,h_n)}{\widehat{f}_{n,0}({ x,h_n})}\cdot
$$ 
It is the local polynomial estimator of order 0 of the conditional cumulative distribution function. 
\item The estimator $\widehat{F}^{(1)}_{n}(t,h_n|{ x})$ is better than the Nadaraya-Watson estimator when the design is random and has the favorable property to reproduce polynomial of order 1.
Precisely, the local linear estimator has a high minimax efficiency among all possible estimators, including nonlinear smoothers (see Fan and Gijbels~\cite{Fan}).
\item We have state in the beginning of this Section that we restrict ourselves to the local polynomial estimator of order 1. 
The local polynomial estimator can be generalized to the orders $p\geqslant 2$, but the equations become more complicated.  
We show briefly the form of the local polynomial estimator of order 2: 
$$
\widehat{F}^{(2)}_{n}(t,h_n|{ x})=\frac{a_1\widehat{r}_{n,0}({ x},t,h_n)+a_2\widehat{r}_{n,1}({ x},t,h_n)+a_3\widehat{r}_{n,2}({ x},t,h_n)}{a_1\widehat{f}_{n,0}({ x,h_n})+a_2\widehat{f}_{n,1}({ x,h_n})+a_3\widehat{f}_{n,2}({ x,h_n})}
$$
where 
$\left\{
\begin{array}{ll}
a_1&= \widehat{f}_{n,2}({ x,h_n})\widehat{f}_{n,4}({ x,h_n})-\left(\widehat{f}_{n,3}({ x,h_n})\right)^2\\ 
 a_2&=\widehat{f}_{n,2}({ x,h_n})\widehat{f}_{n,3}({ x,h_n})-\widehat{f}_{n,1}({x,h_n})\widehat{f}_{n,4}({ x,h_n})\\
 a_3&=\widehat{f}_{n,1}({ x,h_n})\widehat{f}_{n,3}({ x,h_n})-\left(\widehat{f}_{n,2}({ x,h_n})\right)^2
 \end{array}
\right.
$

\noindent 
and $\widehat{f}_{n,3},\widehat{f}_{n,4}$ and $\widehat{r}_{n,2}$ are the direct extensions of the definitions given in the Equations \eqref{defFnj} and \eqref{defRnj}.
Note also that, it is not very interesting to study $p\geqslant 3$, see Fan and Gijbels~\cite{Fan}, pp. 20-22 and 77-80.
The argument is that the mean square error increases with $p$. 
\end{enumerate}
\end{remark}
\vspace{0.3cm}

Now, we study the consistency of the estimator $\widehat F^{(1)}_{n}(t,h_n|x)$ via the following decomposition:  
\[
\widehat F^{(1)}_{n}(t,h_n|x)-F(t|x)=
\underbrace{\widehat F^{(1)}_{n}(t,h_n|x)-\widehat{\mathbb E}\left( \widehat F^{(1)}_{n}(t,h_n|x)\right)}_{(1)}+
\underbrace{\widehat{\mathbb E}\left( \widehat F^{(1)}_{n}(t,h_n|x)\right)-F(t|x)}_{(2)}
\]

\noindent 
where, following the ideas of Deheuvels and Mason (see ~\cite{DeheuvelsMason}), the centering term is defined by:
\[
\widehat{\mathbb{E}}\left(\widehat{F}^{(1)}_{n}(t,h_n|{ x})\right)=\frac{f_{n,2}({ x,h_n})r_{n,0}({ x},t,h_n)-f_{n,1}(x,h_n)r_{n,1}({ x},t,h_n)}{f_{n,0}({ x},h_n)
f_{n,2}({ x},h_n)-f^2_{n,1}({ x},h_n)}
\]
where $f_{n,j}({ x,h_n})=\mathbb{E}\left(\widehat{f}_{n,j}({ x,h_n})\right)$ for $j=0,1,2$ and $r_{n,j}({ x},t,h_n)=\mathbb{E}\left(\widehat{r}_{n,j}({ x},h_n)\right)$ for $j=0,1$.
\vspace{0.3cm}

The  \textit{random part} (1) is the object of our theorem given in the following Section. 
Under (F.1-2), (H.1) and (K.1-3), the {\it deterministic term} (2), so-called bias, converges uniformly to 0 over $(x,t)\in I\times \mathbb R$. The argument to proof this is the Bochner's Lemma (see for instance~\cite{Einmahl}, or our Equations \eqref{E:bochner} in Section \ref{proofs}).

\section{\bf Uniform consistency of the local linear estimator}\label{Sect:Main_Result}
We have now all the ingredients to state our main results. 
The uniform law of the logarithm concerning the local linear estimator of the {\it cond-cdf}, is given in Theorem~\ref{theo1} below. 
\begin{theorem}
\label{theo1}
Under (F.1-2), (H.0-2) and (K.1-3), we have:
\begin{equation}
\sup_{{ x}\in { I}}
\sqrt{\frac{nh_n}{\log(h_n^{-1})}}\left|\widehat{F}^{(1)}_{n}(t,h_n|{ x})-\widehat{\mathbb{E}}\left(\widehat{F}^{(1)}_{n}(t,h_n|{ x})\right)\right| \xrightarrow[n\to+\infty]{\mathbb{P}}\sigma_{F,t}({ I})
\end{equation}
where 
$
\sigma^2_{F,t}({ I})=
2||{ K}||_2^2\sup_{{ x}\in { I}}\left(\frac{F(t|{ x})(1-F(t|{ x}))}{f_{ X}(x)}\right)\cdot
$

Moreover, we have: 
\begin{equation}
\sup_{t\in\mathbb{R}}\sup_{{ x}\in { I}}\sqrt{\frac{nh_n}{\log(h_n^{-1})}}\left|\widehat{F}^{(1)}_{n}(t,h_n|{ x})-\widehat{\mathbb{E}}\left(\widehat{F}^{(1)}_{n}(t,h_n|{ x})\right)\right|\xrightarrow[n\to+\infty]{\mathbb{P}}\sigma_{F}({ I})
\end{equation}
where 
\[
\sigma^2_F({ I})=2||{ K}||_2^2\sup_{t\in\mathbb{R}}\sup_{{ x}\in { I}}\left(\frac{ F(t|{ x})(1-F(t|{ x})) }{f_{ X}(x)}\right)
=\frac{||{ K}||_2^2}{2\displaystyle\inf_{{ x}\in{ I}}f_{ X}({ x})}\cdot
\]
\end{theorem}

\noindent The proof of Theorem~\ref{theo1} is postponed to Section~\ref{proofs}.

\begin{remark}
\hspace{0.1cm}
\begin{enumerate}
\item The matching almost surely convergence result can also be obtained by assuming (H.2-4) instead of (H.0-2). 

\item 
The terms $\sigma_{F,t}({I})$ and $\sigma_{F}({I})$ depend upon the unknown density $f_X$.
But it is a minor problem in practice, because, as shown in Deheuvels~\cite{Deheuvels00}, 
and Deheuvels and Mason~\cite{DeheuvelsMason}, an application of Slutsky's Lemma allows us to replace, without loss of
generality, this quantity by $\widehat f_{n,0}(x,h_n)$ (or by any other estimator of $f_X(x)$ which
is uniformly consistent on $I$). Indeed, under (F.1-2), (H.0-2), (K.1-3) we have $ 
\sup_{x\in I}\left| \frac{\widehat f_{n,0}(x,h_n)}{f_X(x)}-1\right|\xrightarrow[n\to+\infty]{\mathbb{P}}0$. 
\end{enumerate}
\end{remark}

\noindent
This last remark yields to the following corollary.
\begin{corollary}\label{C:bands}
Under (F.1-2), (H.0-2), (K.1-3), we have:
\begin{equation}\label{LoiUnif2}
\sup_{t\in\mathbb{R}}\sup_{{x}\in { I}}\sqrt{\frac{2nh_n}{\|K\|_2^2\log(h_n^{-1})}\widehat f_{n,0}(x,h_n)}\left|\widehat{F}^{(1)}_{n}(t,h_n|{ x})-\widehat{\mathbb{E}}\left(\widehat{F}^{(1)}_{n}(t,h_n|{x})\right)\right|\xrightarrow[n\to+\infty]{\mathbb{P}}1\cdot
\end{equation}
\end{corollary}  
\vspace{0.3cm}

\noindent We introduce the following quantity $L_n(x):=\displaystyle\sqrt{\frac{2nh_n}{\|K\|_2^2\log(h_n^{-1})}\widehat f_{n,0}(x,h_n)}^{\,-1}$.  We have noted at the end of  the Section \ref{Intro} that the bias part can be neglected, then we have the following proposition. 
\begin{proposition}\label{propAurelie}
Under (F.1-2), (H.0-2) and (K.1-3), and if $h_n$ is such that the bias term\\
$\sup_{t\in\R}\sup_{x\in I}\{L_n(x)\}^{-1}\left| F(t|x)-\widehat{\mathbb{E}}\left(\widehat{F}^{(1)}_{n}(t,h_n|{ x})\right)\right|\xrightarrow[n\to+\infty]{}0$ then we have:
\begin{equation}\label{LoiUnif2Bis}
\sup_{t\in\mathbb{R}}\sup_{{ x}\in { I}}\left\{L_n(x)\right\}^{-1}\left|\widehat{F}^{(1)}_{n}(t,h_n|{ x})-F(t|x)\right|\xrightarrow[n\to+\infty]{\mathbb{P}}1.
\end{equation}

\end{proposition}

\begin{remark}
\hspace{0.1cm}
\begin{enumerate}
\item
The matching almost surely convergence result can also be obtained by assuming (H.2-4) instead of (H.0-2). 
\item 

For our applications in Sections \ref{appli} and \ref{Sect:feto}, a reference choice for $h_n$ is given by minimizing the weighted Mean Integrated Square Error (MISE) criteria (see for instance Berlinet~\cite{Berlinet}, Deheuvels~\cite{Deheuvels77} or Deheuvels and Mason~\cite{DeheuvelsMason}). A detailed discussion about the theoretical choice of this bandwidth is given in Ferrigno~\cite{Ferrigno}. The asymptotically optimal constant bandwidth is given by: 
$$h_n=C(K,F,f_X)n^{-\frac{1}{5}}$$
where the constant $C(K,F,f_X)$ is easy to calculate.

%
%
\item 
The choice of the kernel $K$ is not important in practice. 
The most common used kernels are the Gaussian, the indicator function over $[-\frac12, \frac12]$, and the Epanechnikov kernels (see for instance Deheuvels~\cite{Deheuvels77}). Note that the Gaussian kernel is not compactly supported, but our results can be extended to this case. 
\end{enumerate}
\end{remark}
\vspace{0.2cm}

\section{\bf Uniform asymptotic certainty bands}\label{Sect:3}
\subsection{Application to the {\it cond-cdf}}
We show now how the Proposition \ref{propAurelie} can be used to construct uniform asymptotic certainty bands for $F(t| x)$, in the following sense. 
%
%
%
%
%
Under the assumptions of the Proposition 1, we have, for each $0<\varepsilon<1$, and as $n\to+\infty$: 

%
\begin{equation}\label{ConfBands1}
\mathbb P\left\{F(t|x)\in \left[\widehat{F}^{(1)}_{n}(t,h_n|{ x})\pm(1+\varepsilon)L_n(x)\right],\textup{ for all } (x,t)\in I\times \mathbb R\right\}\to 1\end{equation}
 and
\begin{equation}\label{ConfBands2}
\mathbb P\left\{F(t|x)\in \left[\widehat{F}^{(1)}_{n}(t,h_n|{ x})\pm(1-\varepsilon)L_n(x)\right], \textup{ for all }(x,t)\in I\times \mathbb R\right\}\to 0.\end{equation}
 
 \vspace{0.4cm}
%
 
 Whenever (\ref{ConfBands1}) and (\ref{ConfBands2}) hold jointly for each $0<\varepsilon<1$, we have the following corollary:
\begin{corollary}\label{Prop2}
Under (F.1-2), (H.0-2) and (K.1-3), and if $h_n$ is such that the bias term \\
$\sup_{t\in\R}\sup_{x\in I}\{L_n(x)\}^{-1}| F(t,h_n|x)-\widehat{\mathbb{E}}\left(\widehat{F}^{(1)}_{n}(t,h_n|{ x})\right)|\xrightarrow[n\to+\infty]{}0$ then the interval
 \begin{equation}\label{ConfBands4}
\left[\widehat{F}^{(1)}_{n}(t,h_n|{ x})\pm L_n(x)\right]\end{equation}
provides uniform asymptotic certainty bands (at an asymptotic confidence level of 100\%) for the {\it cond-cdf} $F(t|x)$, uniformly in $(x,t)\in I\times \mathbb R.$

\end{corollary}
\vspace{0.3cm}
\begin{remark}
\hspace{0.3cm}
\begin{enumerate}
\item Probability convergence is sufficient for forming certainty bands, and requires less restrictive hypotheses on the bandwidth $h_n$ than the almost surely convergence results. That is why we use only the probability convergence result of the Proposition \ref{propAurelie}. 
\item Following a suggestion of Deheuvels and Derzko~\cite{DeheuvelsDerzko}, we use, for these upper and lower bounds for $F(t|x)$, the qualification of certainty bands, rather that of confidence bands, because there is no preassigned confidence level $\alpha\in (0,1)$. Some authors (see for instance Deheuvels and Mason~\cite{DeheuvelsMason}, or Blondin~\cite{Blondin}) have used the term confidence bands. 
\end{enumerate}
\end{remark}
\subsection{Application to the regression function}
Let $m(x)=\E(Y|X=x)$ the regression function and $\widehat m_{n}^{(1)}(x)=\int y\widehat F_n^{(1)}(dy,h_n|x)$ its local linear estimator. The Proposition \ref{propAurelie} has the following corollary for the regression function.

\begin{corollary}\label{CorRegression}
Under (F.1-2), (H.0-2) and (K.1-3), and if $h_n$ is such that the bias term\\
$\sup_{t\in\R}\sup_{x\in I}\{L_n(x)\}^{-1}\left| F(t|x)-\widehat{\mathbb{E}}\left(\widehat{F}^{(1)}_{n}(t,h_n|{ x})\right)\right|\xrightarrow[n\to+\infty]{}0$ and the variable $Y$ lives in the real interval $[\alpha,\beta]$, then the interval
 \begin{equation}\label{ConfBands6}
\left[\widehat{m}_n^{(1)}( x)\pm (\beta-\alpha)L_n(x)\right]\end{equation}
provides uniform asymptotic certainty bands (at an asymptotic confidence level of 100\%) for the conditional regression function $m(x)$, uniformly in $x\in I.$

\end{corollary}
The proof of Corollary~\ref{CorRegression} is postponed to Section~\ref{proofs}.

\subsection{Application to the conditional quantiles}
Let $0<\alpha<1$. We define the conditional $\alpha$-quantile of the {\it cond-cdf} by:
\[
q_\alpha(x)=\inf\{t\in\R : F(t|x)\geqslant \alpha\},\quad\mbox{for all }\alpha\in (0,1). 
\]
The local linear estimator of the conditional $\alpha$-quantile is defined by: 
\[
\widehat q_{\alpha,n}^{(1)}(x)=\inf\{t\in\R : \widehat F_n^{(1)}(t,h_n|x)\geqslant \alpha\},\quad\mbox{for all }\alpha\in (0,1). 
\]

The Proposition \ref{propAurelie} has the following corollary for the conditional quantiles.

\begin{corollary}\label{CorQuantiles}
Under (F.1-2), (H.0-2) and (K.1-3), if $h_n$ is such that the bias term\\
$\sup_{t\in\R}\sup_{x\in I}\{L_n(x)\}^{-1}\left| F(t|x)-\widehat{\mathbb{E}}\left(\widehat{F}^{(1)}_{n}(t,h_n|{ x})\right)\right|\xrightarrow[n\to+\infty]{}0$ and if the function $x\mapsto f_{X,Y}\left(x,q_\alpha(x)\right)\neq0$ for all $x\in I$, then the interval
 \begin{equation}\label{ConfBands5}
\left[\widehat{q}_{\alpha,n}^{(1)}( x)\pm \frac{2L_n(x)f_X(x)}{f_{X,Y}\left(x,q_\alpha(x)\right)}\right]\end{equation}
provides uniform asymptotic certainty bands (at an asymptotic confidence level of 100\%) for the conditional $\alpha$-quantile $q_\alpha(x)$, uniformly in $x\in I.$

\end{corollary}
The proof of Corollary~\ref{CorQuantiles} is postponed to Section~\ref{proofs}.

\begin{remark}
\hspace{0.3cm}
\begin{enumerate}
\item The form of these certainty bands is not very useful in practice since the bounds depend upon the unknown conditional density $f_{Y|X}(y|x)=\frac{f_{X,Y}(x,y)}{f_X(x)}\cdot$ Nevertheless, this gives the order of the deviation $\left|\widehat q_{\alpha,n}^{(1)}(x)-q(x)\right|$. 
\item To give a more practical result, the idea is to replace the conditional density $f_{Y|X}(q_\alpha(x)|x) $ by an estimator $\widehat f_{Y|X}\left(\widehat q_{\alpha,n}^{(1)}(x)|x\right)$ such that $\sup_{x\in I}\left|\frac{\widehat f_{Y|X}\left(\widehat q_{\alpha,n}^{(1)}(x)|x\right)}{f_{Y|X}(q_\alpha(x)|x)}-1\right|\xrightarrow[n\to+\infty]{\mathbb{P}}0$. This is not the object of the present article, and will be presented in a future work. A review of kernel estimators for the conditional density is given for instance in~\cite{Youndje,Youndje2}. We can cite here the kernel estimator of Parzen-Rosenblatt~\cite{Parzen,Rosenblatt}. 
\end{enumerate}
\end{remark}

\section{{\bf A simulation study}}\label{appli}

In this paragraph, the {\it cond-cdf} and the certainty bands introduced
in Corollary \ref{Prop2} are constructed on simulated data. We considered the
case: $X\sim\mathcal N(0,1)$ where $\mathcal N(0,1)$ denotes  the Gaussian distribution with mean 0 and standard deviation 1. We present two models: 
\begin{enumerate}
\item[(M$_1$)] 
$Y|X=x$ follows a Beta$(a,b)$ distribution with shape parameters $a=1$ and $b=1+x^2$. 
\item[(M$_2$)] 
$Y|X=x$ follows an Uniform distribution between $-|x|$ and $|x|$.
\end{enumerate}

\begin{figure}[htbp]
\begin{center}
\includegraphics[width=0.7\textwidth]{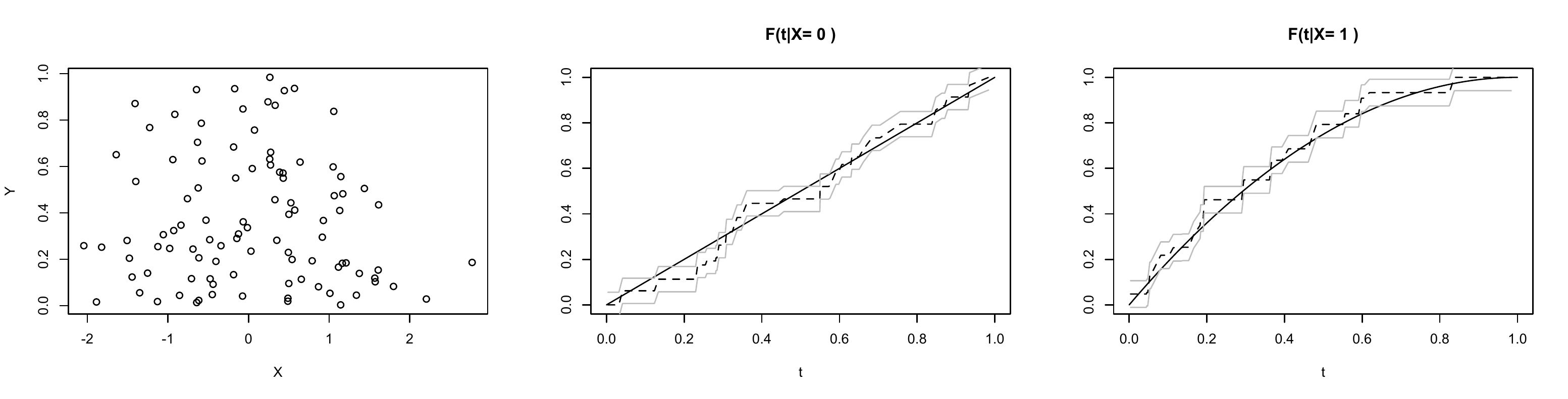}
\includegraphics[width=0.7\textwidth]{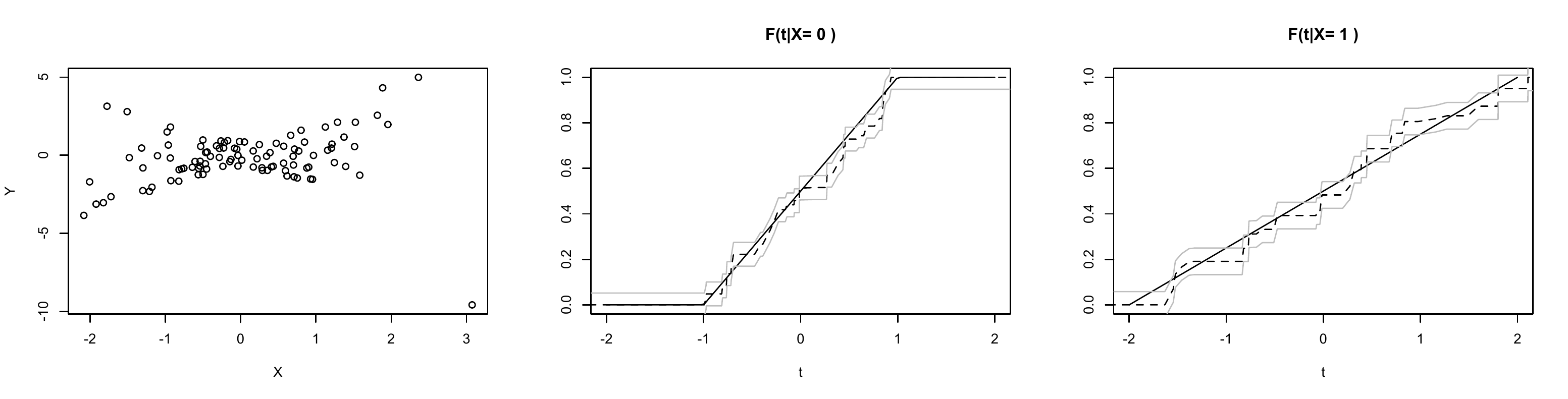}
\includegraphics[width=0.7\textwidth]{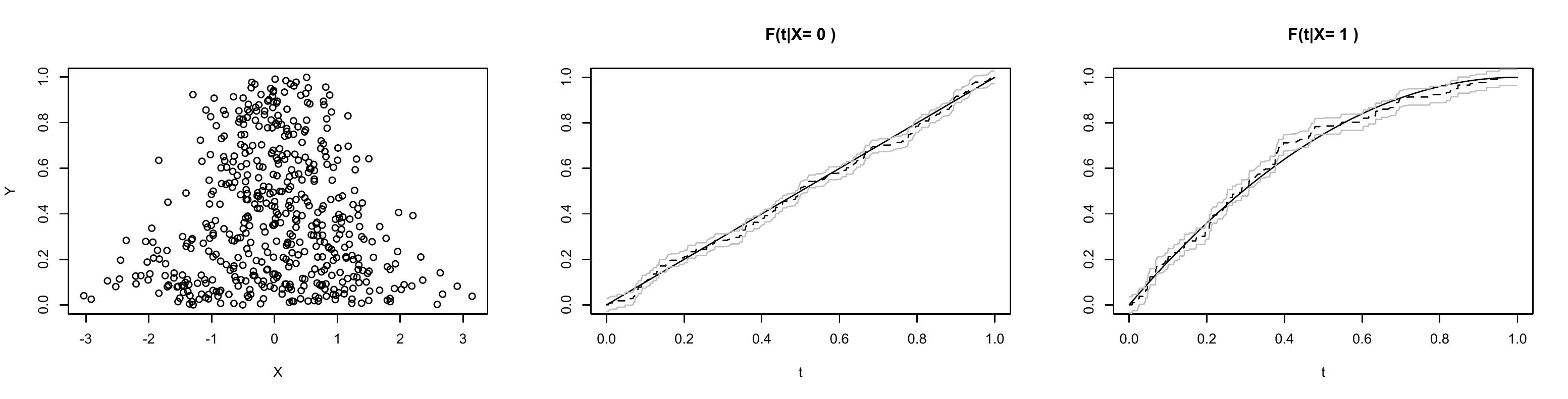}
\includegraphics[width=0.7\textwidth]{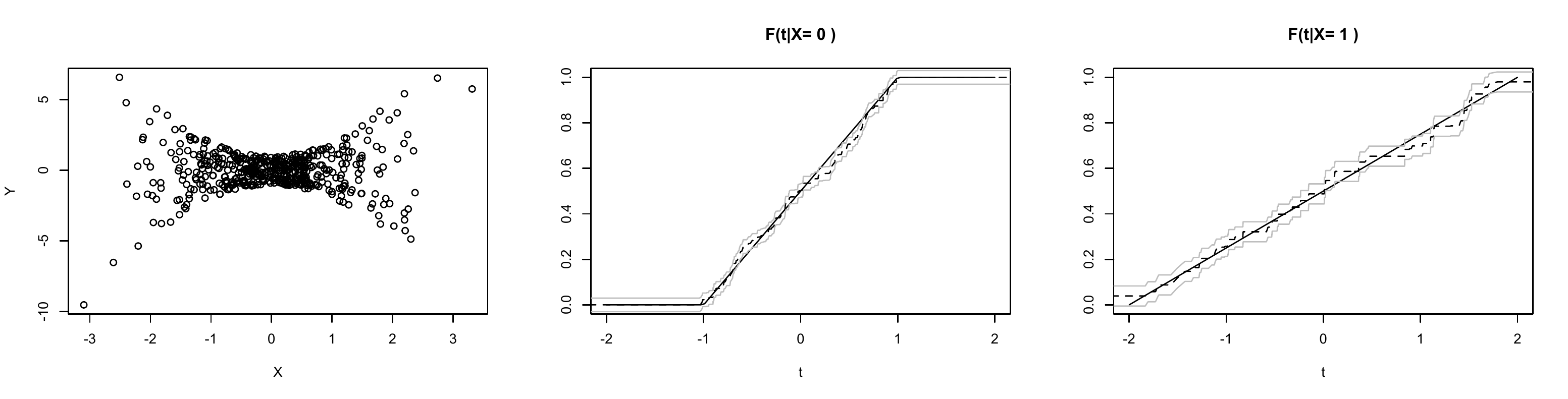}
\caption{From top to bottom: models (M$_1$) and (M$_2$) for $n=100$, and (M$_1$) and (M$_2$) for $n=500$.}
\label{Plot}
\end{center}
\end{figure}

We worked with the sample sizes $n$ = 100 and $n=500$. For the kernel $K$, we opted for the Epanechnikov kernel. For the bandwidth, we selected $h_n =n^{-1/5}$. 
The Figure \ref{Plot} illustrates the results for the models (M$_1$) and (M$_2$) defined above. For each model, we give the graph of a sample $(X_i,Y_i)_{i=1,\ldots,n}$, and the {\it cond-cdf}: the true function  $F(t|x)$ is in full line, whereas the estimated conditional distribution $\widehat F_n^{(1)}(t,h_n|x)$ is in black dashed line, and certainty bands in grey line, for $x=0$ and $1$. 

The confidence bands appear to be adequate. The fact that the true function does not belong to our certainty bands for some points was expected: it is due to the $\varepsilon$ term in Equations \eqref{ConfBands1} and \eqref{ConfBands2}.  For $n=500$, the results are better than for $n=100$. 

\section{{\bf Application in study of the fetal growth}}\label{Sect:feto}
The study is based on 3606 fetuses autopsied in fetopathologic units of the
"Service de foetopathologie et de placentologie" of the Maternité Régionale Universitaire (CHU Nancy, France) between 1996 and 2013. From this dataset, 694 fetuses were carefully selected by exclusion of multiple pregnancies, malformed, macerated or serious ill fetuses, or those with chromosomal abnormalities.

The naive idea, classically used by the fetopathologists or the echographists (see for instance~\cite{Altman},~\cite{Royston}), is to fit a parametric regression model $Y_i=\beta_0+\beta_1  X_i+\beta_2 X_i^2+\epsilon_i$ with the assumptions that $\epsilon_i$, for $i=1,\ldots,n$ are independent and follow the Gaussian distribution $\mathcal N(0,\sigma)$. The parameters $\beta_0,\beta_1,\beta_2,\sigma$ are estimated by the least squares method. We use the \texttt{R 2.15.1} function \texttt{lm}. 

\begin{figure}[htbp]
\begin{center}
\includegraphics[width=0.48\textwidth]{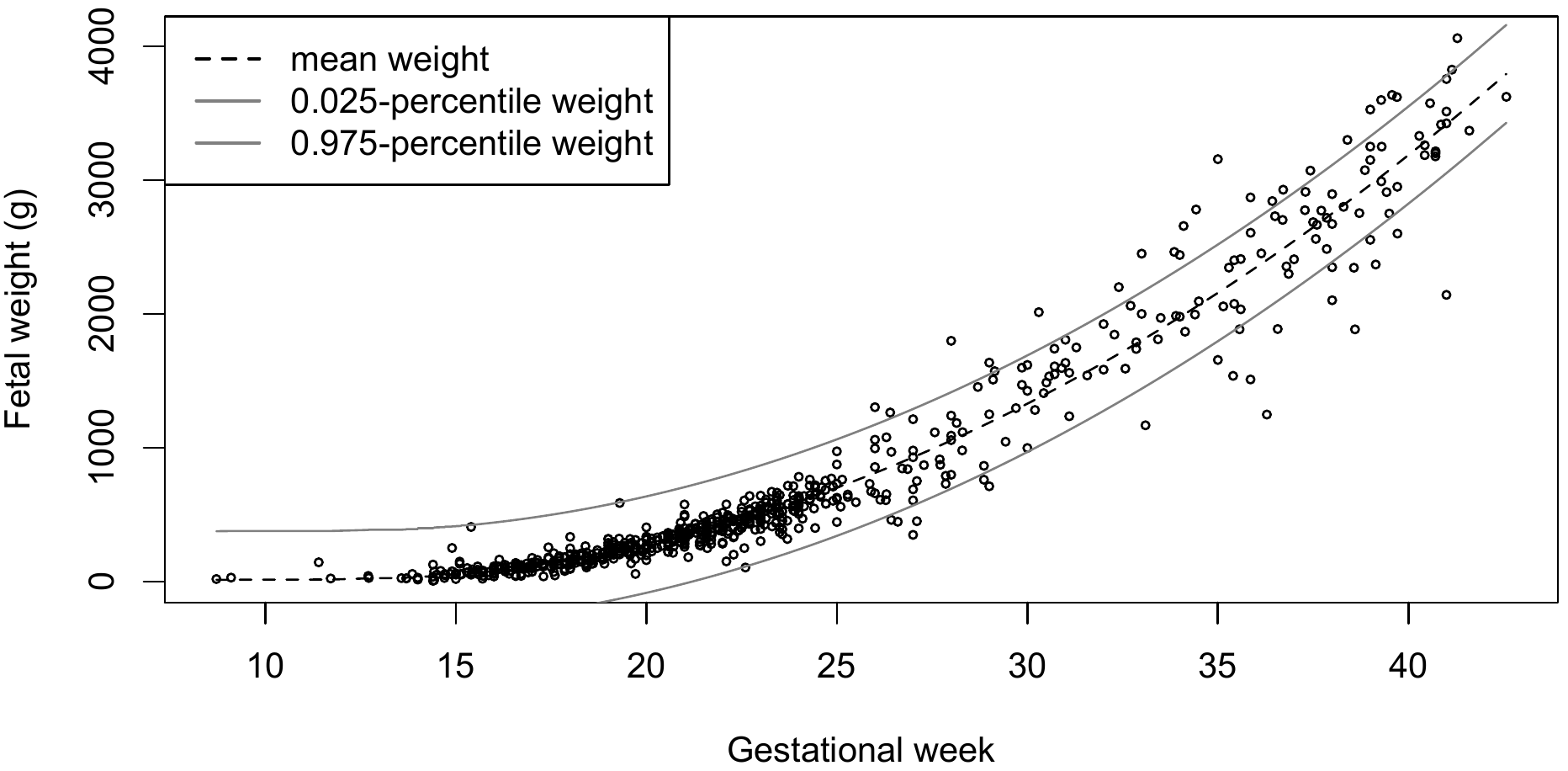}
\includegraphics[width=0.48\textwidth]{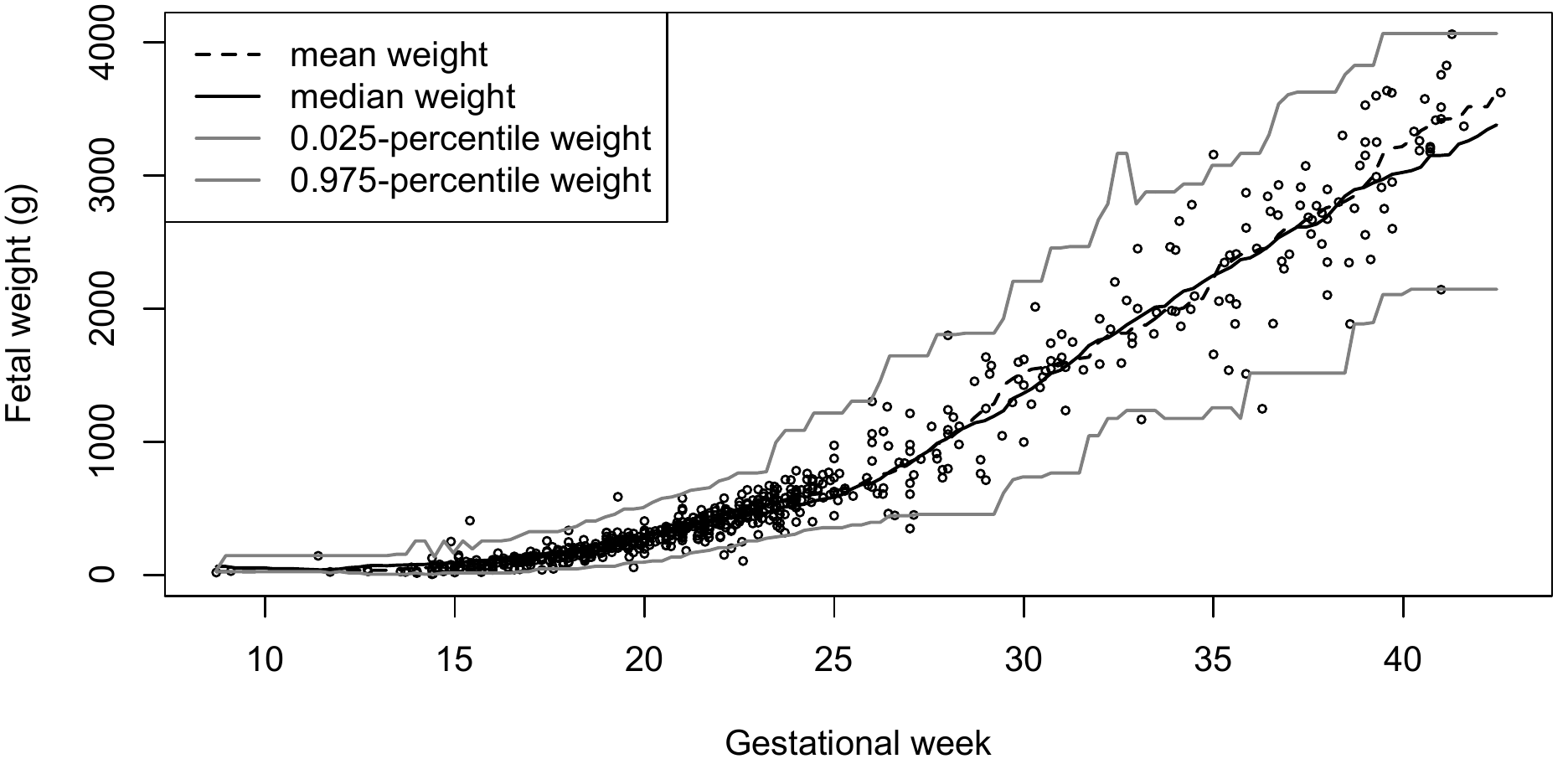}
\caption{Fetal weight during the pregnancy: estimation of mean and quantiles with the second order polynomial regression (left), and with the linear local method (right). }
\label{Fig}
\end{center}
\end{figure}

The result is shown on the left graph of the Figure \ref{Fig}. This method yields to several problems: 
\begin{itemize}
\item We obtain heteroscedastic and non-Gaussian errors. 
\item Moreover, regarding the confidence intervals of the previsions, they show that the prevision uncertainty is not growing with the gestational week: this is not consistent with the medical intuition. 
\item Another problem is that the global polynomial estimation can not enhance some changes in the growing curve of the fetal weight. For the fetopathologists, such changes are important as they correspond to delicate periods during the intrauterine growth. These change points can not been observed by a global estimation. 
\end{itemize}

For these reasons, the local polynomial estimation is then a non-parametric alternative to the global parametric regression model. 

We can conclude, by the observation of the right graph of the Figure \ref{Fig}: 
\begin{itemize}
\item Our method gives the mean, the confidence intervals and the median weight. Satisfactorily, the confidence intervals show the growing of the prevision uncertainty with the gestational week. 
\item We observe for instance a change point between the 20th and 25th gestational week on the 0.975 percentile curve. This change point corresponds to the viability date of the fetus. We can also remark a decrease of the growing speed around the 35th week. This has also been remarked in the medical article~\cite{Guihard}, where it is explained that this time corresponds to the regression (in the medical sense) of the placenta. More precise statistical tests to detect the change points of the fetal growth will be presented in a future work.
 \end{itemize}

\section{{\bf Proofs}}\label{proofs}

 \subsection{Proof of Theorem \ref{theo1}}
We prove  the probability convergence result (Theorem \ref{theo1}), uniformly in $(x,t)\in I\times \mathbb R$. 
The uniform in $x\in I$ result (less difficult) is left to the reader. 

\begin{remark}
The almost surely convergence result could be proved with some additional arguments, based on the Borel-Cantelli's Lemma (see for instance Blondin~\cite{Blondin} or Einmahl and Mason~\cite{Einmahl}). 
\end{remark}

\vspace{0.3cm}

\noindent {\bf Step 1}: In a first step, we introduce a general local empirical process.
For any $j=0,1,2$ and continuous real valued functions $c(\cdot)$ and $d(\cdot)$ on ${ J}$, set for $x\in { J}$, $t \in \mathbb{R}$,
\begin{equation}
W_{n,j}({ x},t)=\sum_{i=1}^n\left(c({ x})\1_{\{Y_i \leqslant t\}}+d({ x})\right)
{ K}_j\left(\frac{{ x}-{ X}_i}{h_n}\right)-n\mathbb{E}\left(\left(c({ x})\1_{\{Y\leqslant t\}}+d({ x})\right)
{ K}_j\left(\frac{{ x}-{ X}}{h_n}\right)\right)
\end{equation}
where 
\begin{equation}\label{E:Kj}
K_j( u)=u^j K( u)
\end{equation} for $j=0,1,2$ and ${ u}\in\mathbb{R}$.

For every fixed $t\in \mathbb{R}$, and $j=0,1,2$, the process $W_{n,j}(\cdot,t)$ can be represented as a bivariate empirical process indexed by a class of functions. More precisely, we have:
$$W_{n,j}(\cdot,t)=\sqrt n\alpha_n(g)=\sum_{i=1}^n \{g(X_i,Y_i)-\mathbb E \left(g(X,Y)\right)\},$$
where $\alpha_n$ is the empirical process based upon $(X_1,Y_1), \ldots, (X_n,Y_n)$ and indexed by a suitable subclass $\mathcal F_{n,j}$ of the class of functions defined on $J\times \mathbb R$:
$$\mathcal F_j =\left\{ (x,y) \mapsto \left\{c(z) \1_{\{y\leqslant t\}} + d(z)\right\} K_j\left(\frac{z-x}{h}\right): t\in \mathbb R, z\in I, 0<h<1\right\}\cdot$$

\vspace{0.3cm}
We give now the result from which our main Theorem \ref{theo1} follows.

\vspace{0.3cm}
\begin{theorem}
\label{theo2}
Under (F.1-2), (H.0-2), (K.1-3) we have:
\begin{equation}
\sqrt{\frac{1}{2nh_n\log(h_n^{-1})}}\sup_{t\in\mathbb{R}}\sup_{x\in I}\left| W_{n,j}(x,t)\right|\xrightarrow[n\to+\infty]{\mathbb{P}}\sigma_{W,j}(I),
\end{equation}
where 
\[
\sigma_{W,j}^2(I)=\sup_{t\in\mathbb{R}}\sup_{x\in I}\mathbb{E}\left(\left|c(x)\1_{\{Y\leqslant t\}}+d(x)\right|^2|X=x\right) f_X(x)\int_{\mathbb{R}}K^2_j(u)du.
\]
\end{theorem}
\vspace{0.3cm}

\noindent 
{\bf Proof}: The proof is divided into an upper bound result, and a lower bound result. 
\vspace{0.3cm}

\textbf{Upper bound part:}
The proof of the upper bound result is divided into two steps. 
The hypothesis (F.2) is important in the upper bound part: we can found a real $M>0$ such that if $X\in J$ then $|Y|\leqslant M$.
\vspace{0.3cm}

\noindent
\textit{Step A: Discretization in $x\in I$ and $t\in [-M,M]$.} 
First, we examine the behavior of our process $(x,t)\mapsto W_{n,j}(x,t)$ on an appropriate chosen grid of $I\times [-M,M]$, with increment $\delta h_n$ for $I$ and increment $\delta$ for $[-M,M]$, for fixed $0<\delta<1$:
$$\left\{\begin{array}{ll} 
 z_{i,n}=a+i\delta h_n,&i=1,\ldots,i_n=\left[\frac{b-a}{\delta h_n}\right],\\
 \\
 t_{l}=-M+l\delta,&  l=1,\ldots, L=\left[\frac{2M}{\delta}\right],\\
\end{array}
\right.$$
where $[u]\leqslant u<[u]+1$ represents the integer part of $u$. 

The study of the supremum on $I\times[-M,M]$ is then reduced to the study of the maximum on a finite number of points. The empirical process is then indexed on the finite class of functions 
$$\mathcal F_{n,j} =\left\{ (x,y) \mapsto \left\{c(z_{i,n}) \1_{y\leqslant t_{l}} + d(z_{i,n})\right\} K_j\left(\frac{z_{i,n}-x}{h_n}\right): i=1,\ldots,i_n, \; l=1,\ldots, L \right\}\cdot$$

The useful tool in this Step is the Bernstein inequality (see for instance Deheuvels and Mason~\cite{DeheuvelsMason}). 

\vspace{0.3cm}
\noindent\textit{Step B: Oscillation.}
Next we study the behavior of our process between the grid points $z_{i,n}$ and $t_{l}$ for $1\leqslant i\leqslant i_n$ and $1\leqslant l\leqslant L$. The objective of this Step is to study the maximal oscillations between the grid points. The useful tool in this Step is the Talagrand inequality for VC Classes (see for instance Talagrand~\cite{Talagrand}, Einmahl and Mason~\cite{Einmahl} or Blondin~\cite{Blondin}). 

\vspace{0.3cm}

\vspace{0.3cm}
\textbf{Lower bound part:}
The lower bound result is proved with technical results based on Poisson processes, and needs the Bochner's Lemma. In this part, the hypothesis (F.1) is particularly important.


\vspace{0.3cm}

\noindent {\bf Step 2}: We give now useful corollaries of Theorem \ref{theo2}. We recall that $K_j$ has been defined in Equation \eqref{E:Kj}, and $||K_j||_{2}^2=\int_\R K_j^2(u)du$. 


\begin{corollary}\label{C:2} Under (F.1-2), (H.0-2) and (K.1-3), we have, by application of Theorem \ref{theo2} with $c(x)=0, d(x)=1, j=0,1,2$:

\begin{equation}\label{E:5}
\sqrt{\frac{n h_n}{2\log(h_n^{-1})}}\sup_{x\in I}| \widehat{f}_{n,j}(x)-f_{n,j}(x)|\xrightarrow[n\to+\infty]{\mathbb{P}}\sigma_{f,j}(I),
\end{equation}

\noindent where 


$$\sigma_{f,j}^2(I)=||K_j||_{2}^{2}\sup_{x\in I}\left\{f_X(x)\right\}\cdot$$

\end{corollary}


\begin{corollary}\label{C:4} Under (F.1-2), (H.0-2) and (K.1-3), we have, by application of Theorem \ref{theo2} with $c(x)=1, d(x)=0, j=0,1$:

\begin{equation}
\sqrt{\frac{n h_n}{2\log(h_n^{-1})}}\sup_{t \in \mathbb{R}}\sup_{x\in I}| \widehat{r}_{n,j}(x,t)-r_{n,j}(x,t)|\xrightarrow[n\to+\infty]{\mathbb{P}}\sigma_{r,j}(I)
\end{equation}

\noindent where 


\begin{eqnarray*}\sigma_{r,j}^2(I)&=&||K_j||_{2}^2\sup_{t \in \mathbb{R}}\sup_{x\in I}\left\{F(t|x) f_X(x)\right\}.\end{eqnarray*}


\end{corollary}



Moreover, under (F.1-2), (H.1) and (K.1-3), the Bochner's Lemma (cf for instance~\cite{Einmahl}) implies, uniformly in $(x,t)\in I\times\mathbb R$:

\begin{equation}\label{E:bochner}\begin{array}{ll}
 f_{n,0}(x)&=f_X(x)\mu_0(K)+ o(1),\\
 f_{n,1}(x)&=f_X(x)\mu_1(K)+ o(1) ,\\
 f_{n,2}(x)&=f_X(x)\mu_2(K)+ o(1) ,\\
 r_{n,0}(x,t)&=f_X(x)F(t|x)+ o(1),\\
 r_{n,1}(x,t)&=f_X(x)F(t|x)\mu_1(K)+ o(1),\end{array}\end{equation}


\noindent where $\mu_j(K)=\int_\mathbb R K_j(u)du$, for $j=0,1,2$. 
The hypotheses (K.2-3) imply that $\mu_0(K)=1$, $\mu_1(K)=0$ and $\mu_2(K)\neq 0$. 
\vspace{0.5cm}

\noindent {\bf Step 3}: In this third step, the deviation $\widehat{F}^{(1)}_{n}(t|{ x})-\widehat{\mathbb{E}}\left(\widehat{F}^{(1)}_{n}(t|{ x})\right)$ can be asymptotically expressed as a linear function of the bivariate empirical process.

\vspace{0.3cm}

\noindent $\widehat{F}^{(1)}_{n}(t|{ x})-\widehat{\mathbb{E}}\left(\widehat{F}^{(1)}_{n}(t|{ x})\right)$

\begin{eqnarray}&=& \frac{\widehat r_{n,0}(x,t)\widehat f_{n,2}(x)-\widehat f_{n,1}(x)\widehat r_{n,1}(x,t)}{\widehat f_{n,2}(x)\widehat f_{n,0}(x)-\widehat {f}_{n,1}^{ 2}(x)}- \frac{ r_{n,0}(x,t)f_{n,2}(x)- f_{n,1}(x)r_{n,1}(x,t)}{ f_{n,2}(x) f_{n,0}(x)- f^2_{n,1}(x)}\\  \nonumber \\
&=& \frac{\widehat r_{n,0}(x,t)\widehat f_{n,2}(x)-r_{n,0}(x,t)f_{n,2}(x)+f_{n,1}(x)r_{n,1}(x,t)-\widehat f_{n,1}(x)\widehat r_{n,1}(x,t)}{\widehat f_{n,2}(x)   \widehat f_{n,0}(x)- \widehat f^2_{n,1}(x)} \label{part1}\\ \nonumber\\
&+& \frac{\left(r_{n,0}(x,t)f_{n,2}(x)\!-\!f_{n,1}(x)r_{n,1}(x,t)\right)\!\!(f_{n,2}(x) f_{n,0}(x)\!-\! f^2_{n,1}(x)\!-\!\widehat f_{n,2}(x)\widehat f_{n,0}(x)\!+\!\widehat f^2_{n,1}(x))}{\left(\widehat f_{n,2}(x)\widehat f_{n,0}(x)-\widehat f^2_{n,1}(x)\right)\left(f_{n,2}(x) f_{n,0}(x)- f^2_{n,1}(x)\right)}\;\;\label{part2}\end{eqnarray}
\vspace{0.3cm}

\noindent 
First, studying the numerator of the expression~\eqref{part1}, we have:
\vspace{0.3cm}

\begin{eqnarray*}\textrm{Num\eqref{part1}}&=&\left(\big(\widehat f_{n,2}(x)-f_{n,2}(x)\big)+ \big(f_{n,2}(x)-f_X(x)\mu_2(K)\big)+f_X(x)\mu_2(K)\right)\left(\widehat r_{n,0}(x,t)-r_{n,0}(x,t)\right)\\ \\
&&+\left(\big(r_{n,0}(x,t)-f_X(x)F(t|x)\big) +f_X(x)F(t|x)\right)\left(\widehat f_{n,2}(x)-f_{n,2}(x)\right)\\ \\
& & + f_{n,1}(x)\left(r_{n,1}(x,t)-\widehat r_{n,1}(x,t)\right)\\ \\
&&+\left(\big(\widehat r_{n,1}(x,t)-r_{n,1}(x,t)\big)+r_{n,1}(x,t)\right)\left(f_{n,1}(x)-\widehat f_{n,1}(x)\right).
\end{eqnarray*}
\vspace{0.3cm}

Let $\displaystyle \beta_n=\sqrt{\frac{n h_n}{2\log(h_n^{-1})}}\cdot$ Because of (H.2), we have $\beta_n\longrightarrow +\infty$ as $n \longrightarrow + \infty.$ Combining with the two previous corollaries and Equations \eqref{E:bochner}, we see that: \begin{equation}\label{Num}\displaystyle\beta_n\sup_{t\in\mathbb R}\sup_{x\in I}\left|\textrm{Num\eqref{part1}}-\widehat g_n(x,t)\right|\xrightarrow[n\to+\infty]{\mathbb{P}}0\end{equation} where $\widehat g_n(x,t)=f_X(x)\mu_2(K)\left(\widehat r_{n,0}(x,t)-r_{n,0}(x,t)\right) -f_X(x)F(t|x)\left(\widehat f_{n,2}(x)-f_{n,2}(x)\right).$

\noindent 
\begin{lemma}\label{lem2}
The denominator of the expression \eqref{part1}, denoted $\textrm{Den\eqref{part1}}$, satisfies:
\begin{equation}\label{Den}
\sup_{x\in I}\left|\textrm{Den\eqref{part1}}-\frac{1}{f_X(x)^2\mu_2(K)}\right|\xrightarrow[n\to+\infty]{\mathbb{P}}0.
\end{equation}
\end{lemma}

\noindent{\bf{Proof}:}
For all $\epsilon>0$, let the event $$\mathcal A_\epsilon = \left\{\sup_{x\in I}\left|\textrm{Den\eqref{part1}}-\frac{1}{f_X(x)^2\mu_2(K)}\right|>\epsilon\right\}$$ and for all $B>0$, let the event $$\mathcal B=\left\{ \sup_{x\in I} \left|\left( \widehat f_{n,2}(x)\widehat f_{n,0}(x)-\widehat f_{n,1}(x)^2\right)f_X(x)^2\mu_2(K)\right|>B\right\}\cdot$$

Then we have:
\[
\mathbb P(\mathcal A_\epsilon)=\mathbb P(\mathcal A_\epsilon \cap \mathcal B)+\mathbb P(\mathcal A_\epsilon \cap \mathcal B^c).
\]
In one hand, we have:
\[
\mathbb P(\mathcal A_\epsilon \cap \mathcal B^c)\leqslant \mathbb P(\mathcal B^c) \leqslant \mathbb P\left(\sup_{x\in I}\left| \widehat f_{n,2}(x)\widehat f_{n,0}(x)-\widehat f_{n,1}(x)^2\right|\leqslant \frac{B}{\displaystyle\inf_{x\in I} f_X(x)^2\mu_2(K)}\right).
\]
Taking $B=\displaystyle\frac12\left(\inf_{x\in I} f_X(x)^2\mu_2(K)\right)^2$, we obtain: 
\[
\mathbb P(\mathcal A_\epsilon \cap \mathcal B^c)\leqslant \mathbb P\left(\sup_{x\in I}\left| \widehat f_{n,2}(x)\widehat f_{n,0}(x)-\widehat f_{n,1}(x)^2-f_X(x)^2\mu_2(K)\right| \geqslant\frac12\inf_{x\in I}f_X(x)^2\mu_2(K)\right)\xrightarrow[n\to+\infty]{}0\cdot
\]
\noindent
This last limit is obtained by the following trivial decomposition: 
\begin{eqnarray*}
&&\widehat f_{n,2}(x)\widehat f_{n,0}(x)-\widehat f_{n,1}(x)^2-f_X(x)^2\mu_2(K)=\\
 &&\left(\big(\widehat f_{n,2}(x)-f_{n,2}(x)\big) +\big(f_{n,2}(x)-f_X(x)\mu_2(K)\big)+f_X(x)\mu_2(K)\right)\\
 &&\times\left(\big(\widehat f_{n,0}(x)-f_{n,0}(x)\big)+\big(f_{n,0}(x)-f_X(x)\big)\right)\\
&&+\left(\big(\widehat f_{n,2}(x)-f_{n,2}(x)\big) +\big(f_{n,2}(x)-f_X(x)\mu_2(K)\big)\right)f_X(x),\\
\end{eqnarray*} 
and by applying the Corollary \ref{C:2} and Equations \eqref{E:bochner}, combined with the boundedness property of $f_X$ on $I$. 
\vspace{0.3cm}

\noindent On the other hand, we have:
$$\mathbb P(\mathcal A_\epsilon \cap \mathcal B)\leqslant \mathbb P\left(\sup_{x\in I}\left| \widehat f_{n,2}(x)\widehat f_{n,0}(x)-\widehat f_{n,1}(x)^2-f_X(x)^2\mu_2(K)\right| >\epsilon B\right)\xrightarrow[n\to+\infty]{}0$$ for the same reason as before. We have then proved \eqref{Den}. $\square$

\vspace{0.3cm}
Combining \eqref{Num} and \eqref{Den}, we have:
\begin{equation}\label{Part1}
\beta_n\sup_{t\in\mathbb R}\sup_{x\in I}\left|\eqref{part1}-\frac{\widehat r_{n,0}(x,t)-r_{n,0}(x,t)}{f_X(x)} -\frac{F(t|x)}{f_X(x)\mu_2(K)}\left(\widehat f_{n,2}(x)-f_{n,2}(x)\right)\right|\xrightarrow[n\to+\infty]{\mathbb{P}}0\cdot
\end{equation}

\vspace{0.3cm}
This last limit is due to the following lemma.

\begin{lemma}\label{lem}
Let $\mathcal I\subset \mathbb R^d$, with $d\in\mathbb N^*$, and $X_n, Z_n, Y_n$ and $Y$ random functions defined on $\mathcal I$ such that
$$\sup_{w\in\mathcal I}\left| X_n(w)-Z_n(w)\right|\xrightarrow[n\to+\infty]{\mathbb{P}}0, \quad \; \sup_{w\in\mathcal I}\left| Y_n(w)-Y(w)\right|\xrightarrow[n\to+\infty]{\mathbb{P}}0,$$  and for enough high $B>0$, 
$\displaystyle\lim_{n\to+\infty}\mathbb P(\sup_{w\in \mathcal I} |X_n(w)|>B)=0$ and $Y$ is bounded on $\mathcal I$. 
\vspace{0.3cm}
Then, 
$$ \displaystyle\sup_{w\in\mathcal I}\left| X_n(w)Y_n(w)-Z_n(w)Y(w)\right|\xrightarrow[n\to+\infty]{\mathbb{P}}0.$$
\end{lemma}

\noindent\textbf{Proof:}
\begin{eqnarray*}
\sup_{\mathcal I} \left| X_nY_n-Z_nY\right|\leqslant \sup_{\mathcal I} \left| X_n\right| \sup_{\mathcal I}\left|Y_n-Y\right| +\sup_{\mathcal I}\left| Y\right|\sup_{\mathcal I}\left|X_n- Z_n\right|,
\end{eqnarray*}

then for all $\epsilon>0,$
\begin{eqnarray*}
\mathbb P\left(\sup_{\mathcal I} \left| X_nY_n-Z_nY\right|>\epsilon\right)\leqslant \mathbb P\left(\sup_{\mathcal I} \left| X_n\right| \sup_{\mathcal I}\left|Y_n-Y\right|>\frac{\epsilon}{2}\right) +\mathbb P\left(\sup_{\mathcal I}\left| Y\right|\sup_{\mathcal I}\left|X_n- Z_n\right|>\frac{\epsilon}{2}\right).
\end{eqnarray*}
Introducing the event $\mathcal B=\displaystyle\left\{\sup_{\mathcal I}|Y|>B\right\} $ with $B>0$, we bound the second term: 
$$\mathbb P\left(\sup_{\mathcal I}\left| Y\right|\sup_{\mathcal I}\left|X_n- Z_n\right|>\frac{\epsilon}{2}\right)\leqslant \mathbb P\left(\mathcal B\right) +\mathbb P\left(\sup_{\mathcal I}|X_n-Z_n|>\frac{\epsilon}{2B}\right)\xrightarrow[n\to+\infty]{}0$$ for enough high $B$.
Now, we bound the first term: 
$$\mathbb P\left(\sup_{\mathcal I}\left| X_n\right|\sup_{\mathcal I}\left|Y_n- Y\right|>\frac{\epsilon}{2}\right)\leqslant \mathbb P\left(\sup_{\mathcal I}|X_n|>B\right) +\mathbb P\left(\sup_{\mathcal I}|Y_n-Y|>\frac{\epsilon}{2B}\right)\xrightarrow[n\to+\infty]{}0
$$
for enough high $B$.
$\square$
\vspace{0.3cm} 

\noindent Let's studying now the second part \eqref{part2} in the expression of $\widehat{F}^{(1)}_{n}(t|{ x})-\widehat{\mathbb{E}}\left(\widehat{F}^{(1)}_{n}(t|{ x})\right)$. Note that the numerator is equal to:
$$\left(r_{n,0}(x,t)f_{n,2}(x)-f_{n,1}(x)r_{n,1}(x,t)\right)\left(f_{n,2}(x) f_{n,0}(x)- f^2_{n,1}(x)-\widehat f_{n,2}(x)\widehat f_{n,0}(x)+\widehat f^2_{n,1}(x)\right)\cdot$$

\noindent 
The first term $\textrm{Num}_1\eqref{part2} =r_{n,0}(x,t)f_{n,2}(x)-f_{n,1}(x)r_{n,1}(x,t)$ of this last expression converges uniformly, thanks to \eqref{E:bochner}: 
\begin{equation}\label{Num2_a}
\sup_{t\in\mathbb R}\sup_{x\in I}\left|\textrm{Num}_1\eqref{part2} -f_X(x)^2\mu_2(K)F(t|x)\right|\xrightarrow[n\to+\infty]{}0.
\end{equation}

\noindent With the same arguments as in the study of the numerator of \eqref{part1}, we study the second term defined by:
$$\textrm{Num}_2\eqref{part2} =f_{n,2}(x) f_{n,0}(x)- f^2_{n,1}(x)-\widehat f_{n,2}(x)\widehat f_{n,0}(x)+\widehat f^2_{n,1}(x)$$ and show that

\begin{equation}\label{Num2_b}
\beta_n\sup_{x\in I} \left|\textrm{Num}_2\eqref{part2} -f_X(x)\mu_2(K)\left(f_{n,0}(x)-\widehat f_{n,0}(x)\right) -f_X(x)\left(f_{n,2}(x)-\widehat f_{n,2}(x)\right)\right| \xrightarrow[n\to+\infty]{\mathbb P}0\cdot
\end{equation}

Thanks to the  boundedness property of $f_X$ on $I$, and the Corollary \ref{C:2}, we have, by the Lemma \ref{lem}: 
$$
\beta_n\sup_{t\in\mathbb R}\sup_{x\in I} \left|\textrm{Num}\eqref{part2} -\widehat j_n(x,t)\right| \xrightarrow[n\to+\infty]{\mathbb P}0$$ where $ \widehat j_n(x,t)=f_X(x)^3\mu_2(K)F(t|x)\left(\mu_2(K)\left(f_{n,0}(x)-\widehat f_{n,0}(x)\right)+\left(f_{n,2}(x)-\widehat f_{n,2}(x)\right)\right).
$

\vspace{0.3cm}

\noindent The denominator in~\eqref{part2} can be expressed as
$$\textrm{Den}\eqref{part2} =\textrm{Den}\eqref{part1}\frac{1}{f_{n,2}(x)f_{n,0}(x)-f_{n,1}(x)^2}\cdot$$

It is clear that, thanks to \eqref{E:bochner} 
$$\sup_{x \in I}\left|\frac{1}{f_{n,2}(x)f_{n,0}(x)-f_{n,1}(x)^2}-\frac{1}{f_X(x)^2\mu_2(K)}\right| \xrightarrow[n\to+\infty]{}0\cdot$$

Then, tanks to the boundedness property of $f_X$ on $I$, the Lemma $\ref{lem}$ says that 
\begin{equation}\label{DenPart2}
\sup_{x \in I}\left|\textrm{Den}\eqref{part2}-\frac{1}{f_X(x)^4\mu_2(K)^2}\right| \xrightarrow[n\to+\infty]{\mathbb P}0\cdot
\end{equation}

Finally, we have, thanks to the Lemma \ref{lem}

\begin{equation}\label{Part2}
\beta_n\sup_{t\in\mathbb R}\sup_{x \in I}\left|\eqref{part2}-\frac{F(t|x)}{f_X(x)}\left(f_{n,0}(x)-\widehat f_{n,0}(x)\right) -\frac{F(t|x)}{f_X(x)\mu_2(K)}\left(f_{n,2}(x)-\widehat f_{n,2}(x)\right)\right| \xrightarrow[n\to+\infty]{\mathbb P}0\cdot
\end{equation}

Remember that $\widehat{F}^{(1)}_{n}(t|{ x})-\widehat{\mathbb{E}}\left(\widehat{F}^{(1)}_{n}(t|{ x})\right)=\eqref{part1}+\eqref{part2}$. Then, combining \eqref{Part1} and \eqref{Part2}  it is easy to see finally that 
\begin{equation}\label{E:final}
\beta_n\sup_{t\in\mathbb R}\sup_{x\in I}\left|\widehat{F}^{(1)}_{n}(t|{ x})-\widehat{\mathbb{E}}\left(\widehat{F}^{(1)}_{n}(t|{ x})\right) +\frac{F(t|x)}{f_X(x)}\left(\widehat f_{n,0}(x)- f_{n,0}(x)\right) -\frac{\widehat r_{n,0}(t,x)-r_{n,0}(t,x)}{f_X(x)}\right|  \xrightarrow[n\to+\infty]{\mathbb P}0\cdot
\end{equation}

\vspace{0.3cm}

\noindent {\bf Step 4}:

\vspace{0.3cm}

\noindent Now choosing $\displaystyle c(x)=\frac{1}{f_{X}(x)}$ and $\displaystyle d(x)=-\frac{F(t|x)}{f_{X}(x)}$ in the definition of $W_{n,0}(x,t)$, the local empirical process, it is easy to show that:

\begin{eqnarray*}\label{E:6}W_{n,0}(x,t)&=& \frac{nh}{f_{X}(x)}\widehat r_{n,0}(x,t)-\frac{F(t|x)}{f_{X}(x)} \times nh \widehat f_{n,0}(x)-\frac{nh}{f_{X}(x)} r_{n,0}(x,t)+\frac{F(t|x)}{f_{X}(x)} \times nh f_{n,0}(x)\\ \\
&=& \frac{nh}{f_{X}(x)}\left(\widehat r_{n,0}(x,t)-r_{n,0}(x,t)-F(t|x)\left(\widehat f_{n,0}(x)-f_{n,0}(x)\right)\right)\cdot
\end{eqnarray*}

\noindent Let $A_n=\sqrt{2nh_n \log (h_n^{-1})}^{-1}$, so we have:
\begin{eqnarray}
A_n W_{n,0}(x,t)=\beta_n\left(\frac{\widehat r_{n,0}(t,x)-r_{n,0}(t,x)}{f_X(x)}-\frac{F(t|x)}{f_X(x)}\left(\widehat f_{n,0}(x)- f_{n,0}(x)\right) \right).
\end{eqnarray}

\vspace{0.3cm}

\noindent Then 
\begin{eqnarray*}
&&\beta_n \left(\widehat{F}^{(1)}_{n}(t|{ x})-\widehat{\mathbb{E}}\left(\widehat{F}^{(1)}_{n}(t|{ x})\right)\right)=A_nW_{n,0}(x,t)\\
&&+\beta_n\left(\widehat{F}^{(1)}_{n}(t|{ x})-\widehat{\mathbb{E}}\left(\widehat{F}^{(1)}_{n}(t|{ x})\right) +\frac{F(t|x)}{f_X(x)}\left(\widehat f_{n,0}(x)- f_{n,0}(x)\right) -\frac{\widehat r_{n,0}(t,x)-r_{n,0}(t,x)}{f_X(x)}\right),
\end{eqnarray*}

\noindent and applying Theorem \ref{theo2} 
\begin{equation}
\sup_{t\in\mathbb R}\sup_{x\in I}\left|\beta_n \left(\widehat{F}^{(1)}_{n}(t|{ x})-\widehat{\mathbb{E}}\left(\widehat{F}^{(1)}_{n}(t|{ x})\right)\right)\right| \xrightarrow[n\to+\infty]{\mathbb P}\sigma_F(I)
\end{equation}
\vspace{0.3cm}

\noindent with

\begin{eqnarray*}\sigma^2_{F}(I)&=& \sup_{t \in \mathbb{R}}\sup_{x\in I}\mathbb{E}\left(\left[\frac{\1_{\{{ Y}\leqslant { t}\}}}{f_{X}(x)}-\frac{F(t|x)}{f_{X}(x)}\right]^2|X=x\right)f_{X}(x)||K||^2_2\\ \\
&=& \sup_{t \in \mathbb{R}}\sup_{x\in I}\frac{\mathbb{E}\left(\left[\1_{\{{ Y}\leqslant { t}\}}   - \mathbb{E}({ Y}\leqslant { t}|X=x)\right]^2|X=x\right)}{f_{X}(x)}||K||^2_2\\ \\
&=& \sup_{t \in \mathbb{R}}\sup_{x\in I}\left(\frac{F(t|x)(1-F(t|x))}{f_{X}(x)}\right)||K||^2_2.
\end{eqnarray*}

\noindent This finishes the proof of Theorem \ref{theo1}.$\square$

\subsection{Proof of Corollary \ref{CorRegression}}

 If the variable $Y $ takes its values in $[\alpha,\beta]$, then
$$
\widehat m_n^{(1)}(x)-m(x)=\displaystyle\int_\alpha^\beta y \left(\widehat F^{(1)}_n(dy,h_n|x)-F(dy|x)\right)=-\displaystyle\int_\alpha^\beta\left(\widehat F^{(1)}_n(y,h_n|x)-F(y|x)\right)dy.$$ 

It implies $\left| \widehat m^{(1)}_n(x)-m(x)\right|\leqslant \left| \beta-\alpha\right| \sup_{y\in \R}\left|\widehat F^{(1)}_n(y,h_n|x)-F(y|x)\right|$. 

The conclusion of the Corollary \ref{CorRegression} is then deduced from this last inequality combined with Corollary \ref{Prop2}.

\subsection{Proof of Corollary \ref{CorQuantiles}}

Let $\alpha\in (0,1)$. It suffices to show that \begin{equation}\label{E:CorQuantiles}
\mathbb P\left(\forall x\in I, \; \left| \widehat q^{(1)}_{\alpha,n}(x)-q_\alpha(x)\right|>\frac{2L_n(x)f_X(x)}{f_{X,Y}(x,q_\alpha(x))}\right)\longrightarrow 0\textup{ as }n \longrightarrow + \infty.\end{equation}

Let $\epsilon=\epsilon_{\alpha,x,n}=\frac{2L_n(x)f_X(x)}{f_{X,Y}(x,q_\alpha(x))}\cdot$ We have $$\begin{array}{lcl}
\mathbb P\left(\forall x\in I, \;\left|  \widehat q^{(1)}_{\alpha,n}(x)-q_\alpha(x)\right|>\epsilon\right)& =&\mathbb P\left(\forall x\in I, \;  \widehat q^{(1)}_{\alpha,n}(x)>q_\alpha(x)+\epsilon\right) +\mathbb P\left(\forall x\in I, \; \widehat q^{(1)}_{\alpha,n}(x)<q_\alpha(x)-\epsilon\right)\\
&=&(I)+(II).
\end{array}$$

To study the first term $(I)$, consider the following implication, for all $x\in I$,  $$\left(\widehat q_{\alpha,n}^{(1)}(x)=\inf\{t\in\R : \widehat F_n^{(1)}(t,h_n|x)\geqslant \alpha\}>q_\alpha(x)+\epsilon\right) \Rightarrow\left(\widehat F^{(1)}_{n}(q_\alpha(x)+\epsilon, h_n|x)<\alpha\right).$$

Then 
$$\begin{array}{lcl}
(I)&\leqslant & \mathbb P\left(\forall x\in I,\; \widehat F^{(1)}_{n}(q_\alpha(x)+\epsilon, h_n|x)<\alpha\right)\\
&=& \mathbb P\left(\forall x\in I,\; \widehat F^{(1)}_{n}(q_\alpha(x)+\epsilon, h_n|x)-F(q_\alpha(x)+\epsilon|x)<\alpha-F(q_\alpha(x)+\epsilon|x)\right)\\
&\leqslant&\mathbb P\left(\forall x\in I,\; \sup_{y\in\R}\left|\widehat F^{(1)}_{n}(y, h_n|x)-F(y|x)\right|>F(q_\alpha(x)+\epsilon|x)-\alpha\right)\\
&\leqslant & \mathbb P\left(\forall x\in I,\; \sup_{y\in\R}\left|\widehat F^{(1)}_{n}(y, h_n|x)-F(y|x)\right|>\frac{\epsilon f_{X,Y}(x,q_\alpha(x))}{2f_X(x)}\right)\\
&= &  \mathbb P\left(\forall x\in I,\; \sup_{y\in\R}\left|\widehat F^{(1)}_{n}(y, h_n|x)-F(y|x)\right|>L_n(x)\right)\longrightarrow0\textup{ as } n\longrightarrow +\infty
\end{array}$$
where $F(q_\alpha(x)+\epsilon|x)-\alpha=F(q_\alpha(x)+\epsilon|x)-F(q_\alpha(x)|x)=\epsilon \frac{f_{X,Y}(x,q_\alpha(x))}{f_X(x)}+o(\epsilon)\geqslant\frac{\epsilon f_{X,Y}(x,q_\alpha(x))}{2f_X(x)}\cdot$

To study the second term $(II)$, consider the following implications, for all $x\in I$:
$$\begin{array}{lcl}
\left( \widehat q_{\alpha,n}^{(1)}(x)<q_\alpha(x)-\epsilon \right)&\Rightarrow& \left(F\left(\widehat q_{\alpha,n}^{(1)}(x)|x\right) < F(q_\alpha(x)-\epsilon|x)<\alpha\leqslant \widehat F_n^{(1)}\left(\widehat q_{\alpha,n}^{(1)}(x),h_n|x\right)\right)\\
&&\textup{by the growing property of $F$ and the definition of the quantile $\widehat q_{\alpha,n}^{(1)}(x)$}\\
\\
&\Rightarrow&\left(\widehat F_n^{(1)}\left(\widehat q_{\alpha,n}^{(1)}(x),h_n|x\right) -F\left(\widehat q_{\alpha,n}^{(1)}(x)|x\right)\geqslant \alpha- F(q_\alpha(x)-\epsilon|x)\right).\end{array} $$

Then 
$$\begin{array}{lcl}
(II)&\leqslant &\mathbb P\left(\forall x\in I,\; \widehat F_n^{(1)}\left(\widehat q_{\alpha,n}^{(1)}(x),h_n|x\right) -F\left(\widehat q_{\alpha,n}^{(1)}(x)|x\right)\geqslant \alpha- F(q_\alpha(x)-\epsilon|x)\right)\\
&\leqslant&\mathbb P\left(\forall x\in I,\; \sup_{y\in\R}\left|\widehat F^{(1)}_{n}(y, h_n|x)-F(y|x)\right|>\alpha-F(q_\alpha(x)-\epsilon|x)\right)\\
&\leqslant & \mathbb P\left(\forall x\in I,\; \sup_{y\in\R}\left|\widehat F^{(1)}_{n}(y, h_n|x)-F(y|x)\right|>\frac{\epsilon f_{X,Y}(x,q_\alpha(x))}{2f_X(x)}\right)\\
&= &  \mathbb P\left(\forall x\in I,\; \sup_{y\in\R}\left|\widehat F^{(1)}_{n}(y, h_n|x)-F(y|x)\right|>L_n(x)\right)\longrightarrow0\textup{ as } n\longrightarrow +\infty
\end{array}$$
where $\alpha-F(q_\alpha(x)-\epsilon|x)=F(q_\alpha(x)|x)-F(q_\alpha(x)-\epsilon|x)=\epsilon \frac{f_{X,Y}(x,q_\alpha(x))}{f_X(x)}+o(\epsilon)\geqslant\frac{\epsilon f_{X,Y}(x,q_\alpha(x))}{2f_X(x)}\cdot$

Finally, we have proved \eqref{E:CorQuantiles}. 

\section*{Acknowledgement}
We thank to Doctors Jean-Pierre Masutti and Alain Miton of the Service de foetopathologie et de placentologie of the Maternité Régionale Universitaire (CHU Nancy, France), for the fetal data.


\begin{thebibliography}{00}
\bibitem{Altman}
D.G. Altman, L.S. Chitty, {\it Charts of fetal size: 1. Methodology.}, Br. J. Obstet. Gynaecol., 101 (1994), pp. 29-34.\\
\bibitem{Berlinet}
A. Berlinet, L. Devroye, {\it A comparison of kernel density estimates}, Publ. Inst. Statist. Univ. Paris, 38 (1994), pp. 3-79. \\
\bibitem{Blondin}
D. Blondin, {\it Lois limites uniformes et estimation non paramétrique de la régression}, PhD thesis, Universit\'e de Paris 6, France, 2004.\\
\bibitem{Collomb80}
G. Collomb, {\it Estimation non paramétrique de probabilités conditionnelles}, C. R. Acad. Sci. Paris S\'er. A-B, 291 (1980), pp. 427-430.\\
\bibitem{Collomb}
G. Collomb, {\it Estimation non-paramétrique de la régression: revue bibliographique}, Int. Statist. Rev., 49 (1981), pp. 75-93.\\

\bibitem{Deheuvels77}
P. Deheuvels, {\it Estimation non-paramétrique de la densité par histogramme généralisé}, La Revue de Statistique Appliquée, 35 (1977), pp. 5-42.\\

\bibitem{Deheuvels00}
P. Deheuvels, {\it Limit laws for kernel density estimators for kernels with unbounded supports}, In Asymptotics in Statistics and Probability,
(Ed., M. L. Puri), V.S.P., Amsterdam, 2000.\\

\bibitem{DeheuvelsDerzko}
P. Deheuvels, G. Derzko, {\it Asymptotic certainty bands for kernel density estimators based upon a bootstrap resampling scheme}, Statistical Models and Methods for Biomedical and Technical Systems, 3 (2008), pp. 171-186.\\
\bibitem{DeheuvelsMason}
P. Deheuvels, D.M. Mason, {\it General asymptotic confidence bands based on kernel-type function estimators}, Statist. Infer. Stochastic Process, 7(3) (2004), pp. 225-277.\\
\bibitem{Devroye}
L. Devroye, {\it The uniform convergence of the Nadaraya-Watson regression function estimate}, Can. J. Statist., 6 (1978), pp. 179-191.\\
\bibitem{DonyEinmahlMason}
J. Dony, U. Einmahl, D.M. Mason, {\it Uniform in Bandwidth Consistency of Local Polynomial
Regression Function Estimators}, Austrian Journal of Statistics, 35(2-3) (2006), pp.105-120.\\
\bibitem{Einmahl}
U. Einmahl, D.M. Mason, {\it An empirical process approach to the uniform consistency of kernel-type function estimators}, J. Theoret. Probab., 13(1) (2000), pp. 1-37.\\
\bibitem{Einmahl05}
U. Einmahl, D.M. Mason, {\it Uniform in bandwidth consistency of kernel-type function estimators}, Ann. Statist., 33(3)  (2005), pp. 1380-1403.\\
\bibitem{Fan}
J. Fan, I. Gijbels, {\it Local polynomial modeling and its applications}. Monographs on Statistics and Applied Probability, Chapman and Hall, Vol. 66, 1996.\\ 
\bibitem{Ferrigno}
S. Ferrigno, {\it Un test d'adéquation global pour la fonction de répartition conditionnelle}, PhD thesis, Universit\'e de Montpellier 2, France, 2004.\\
\bibitem{Gannoun}
A. Gannoun, S. Girard, C. Guinot, J. Saracco, {\it References curves based on nonparametric quantile regression}, Statistics in Medecine, 21 (2002), pp. 3119-3155.\\
\bibitem{Guihard}
A.M. Guihard-Costa, {\it Les variations des vitesses de croissance au cours de la vie foetale}, Bulletins et Mémoires de la
Société d'anthropologie de Paris, Nouvelle Série, 5(1-2) (1993), pp.11-20.\\
\bibitem{Hardle}
W. H\"ardle, {\it Applied Nonparametric Regression}, Cambridge University Press, Cambridge, 1990.\\
\bibitem{Hardle88}
W. H\"ardle, P. Janssen, R. Serfling, {\it Strong uniform consistency rates for estimators of conditional functionals}, Ann. Statist., 16(4) (1988), pp.1428-1449.\\
\bibitem{Li}
Q. Li, J. Lin, J. S. Racine, {\it Optimal Bandwidth Selection for Nonparametric Conditional Distribution and Quantile Functions}, Journal of Business and Economic Statistics,
31(1) (2013), pp.57-65.\\
\bibitem{Mint El Mouvid}
M. Mint El Mouvid, {\it Sur l'estimateur lin\'eaire local de la fonction de r\'epartition conditionnelle}, PhD thesis, Universit\'e de Montpellier 2, France, 2000.\\
\bibitem{Nadaraya}
E.A. Nadaraya, {\it On estimating regression}, Theor. Probab. Appl., 9 (1964), pp. 141-142.\\
\bibitem{Parzen}
E. Parzen, {\it On estimation of a probability density function and mode}, Ann. Math. Statist., 33 (1962), pp. 1065-1076.\\
\bibitem{Rosenblatt}
M. Rosenblatt, {\it Remarks on some nonparametric estimates of a density function}, Ann. Math. Statist., 27 (1956), pp. 832-837.\\

\bibitem{Royston}
P. Royston, E.M. Wright, {\it How to construct "normal ranges" for fetal variables}, Ultrasound Obstet. Gynecol., 11 (1998), pp. 30-8.\\

\bibitem{Stute82}
W. Stute, {\it A law of the iterated logarithm for kernel density estimators}, Ann. Probab., 10(2) (1982), pp. 414-422.\\
\bibitem{Stute}
W. Stute, {\it On almost sure convergence of conditional empirical distribution functions}, Ann. Probab., 14(3) (1986), pp. 891-901.\\
\bibitem{Talagrand}
M. Talagrand, {\it Sharper bounds for Gaussian and empirical processes}, Ann. Probab., 22 (1994), pp. 28-76.\\
\bibitem{Tsybakov}
A.B. Tsybakov, {\it Introduction to Nonparametric Estimation}, Springer Series in Statistics, 2008.\\
\bibitem{Watson}
G.S. Watson, {\it Smooth regression analysis}. Sankhy$\overline{a}$ Ser.A, 26 (1964), pp. 359-372.\\
\bibitem{Youndje2}
E. Youndje, {\it Convergence properties of the kernel estimator of conditional density}, Rev. Roumaine Math. Pures Appl., 41(7-8), 1996.\\
\bibitem{Youndje}
E. Youndje, {\it Contribution à l'estimation non-paramétrique par la
méthode du noyau}, Habilitation à diriger des recherches, Université de Rouen et du Havre, France,  2011.\\
\end{thebibliography}
\end{document}